\documentclass{article}
\usepackage{graphicx} 
\usepackage{multicol,amsmath,graphics, mathtools, cancel,soul,color,bbm,amsfonts,dsfont,tikz,mathrsfs,amssymb,bm,hyperref, amsthm}
\usepackage[left=1in, right=1in, top=1.1in,bottom=1.1in]{geometry}

\hypersetup{
    colorlinks=true, 
    linktoc=all,     
    linkcolor=blue,  
}

\usepackage{graphicx}
\usepackage{subcaption}

\newtheorem{theorem}{Theorem}[section]
\newtheorem{corollary}[theorem]{Corollary}
\newtheorem{proposition}[theorem]{Proposition}
\newtheorem{lemma}[theorem]{Lemma}

\title{The coalescent of a sample from a linear-fractional branching process}
\author{Natalia Cardona-Tobón\footnote{Departamento de Estadística, Universidad Nacional de Colombia. C.P. 16486, Bogotá D.C., Colombia. E-mail: {\tt ncardonat@unal.edu.co}}\quad  and \quad  Sandra Palau\footnote{IIMAS, Universidad Nacional Autónoma de México. C.P. 04510, Ciudad de México, México. E-mail: {\tt sandra@sigma.iimas.unam.mx}}}

\begin{document}
\maketitle

\begin{abstract}
In this article, we focus on Bienaymé-Galton-Watson processes with linear-fractional offspring distributions. At a fixed generation, we consider a sample of the individuals alive, drawn in two different ways: either through Bernoulli sampling, where each individual is selected independently with a given probability, or through uniform sampling, where a fixed number of individuals are chosen uniformly at random. We analyze the genealogical trees generated by the sampled individuals. In particular, we establish a relationship between the distributions of the trees resulting from the two sampling schemes.\\


\textbf{MSC2020 subject classifications:} 
60J80; 
60J85; 
62M20; 
92D25. 

\end{abstract}

\section{Introduction}

Some species in nature have cyclical life histories in which reproduction is limited to short, highly synchronized periods separated by prolonged phases of developmental stagnation. An interesting example is the periodical cicadas of the eastern United States (see, e.g., \cite{Fujisawa2018}). These insects are notable for having a long juvenile period that lasts 13 or 17 years depending on the species, during which time they live underground, feed, and remain reproductively inactive. At the end of this juvenile period, all individuals in the given cohort emerge simultaneously in a mass emergence event. 
Adults live on the surface for only a few weeks, during which time they mate, lay eggs, and die. The eggs hatch, and the new generation is buried to begin the next 13- or 17-year cycle.

Branching processes would seem to be a natural framework for modeling the family tree formed by the female descendants of an initial female cicada. In particular, the discrete-time Bienaymé-Galton-Watson (for short BGW) process provides a suitable model for capturing generational reproduction and lineage dynamics. Compared to continuous-time processes, such as birth-death models or diffusions, BGW processes more closely reflect the structured timing of reproduction observed in periodical cicadas, whose generations are clearly separated by long intervals of non-reproductive development.

To analyze the genealogical structure in this setting, we focus on a particular family of offspring distributions known as linear fractional distributions which are mixtures of geometric and Bernoulli distributions. BGW processes with these offspring distributions possess an interesting mathematical property very important in our results, they induce a Markovian structure in the associated coalescent point process (see Section~\ref{sec:CPP}). In addition, we believe that linear-fractional distributions may also be of interest in biology, as they have been used in contexts such as modeling the demography of a rare allele (see, for example, \cite{rannala1997genealogy}, \cite{thompson1976estimation}).

In this paper, we consider a population whose size evolves according to a \textit{Bienaymé-Galton-Watson process} $\{Z_n, n\geq 0\}$, where the offspring distribution is determined by a random variable $\xi$ following a $(p,r)$-linear fractional distribution. In other words, 
\begin{equation}
    Z_0=1 \qquad \text{and}\qquad Z_n= \sum_{i=1}^{Z_{n-1}}\xi_i^{(n)}, \qquad n\ge 1,
\end{equation}
where $\{\xi_i^{(n)}, i, n\ge 1\}$ is a sequence of independent, identically distributed random variables with the same distribution as $\xi$, i.e.
\begin{equation}\label{eq:fractional}
 \mathbb{P}(\xi =0)=1-r \qquad \mbox{and} \qquad \mathbb{P}(\xi =k) =rp(1-p)^{k-1}, \qquad k\geq 1,   
\end{equation}
where $0\leq r \leq 1$ and $0<p<1$. 
We assume that the BGW process is supercritical, i.e.  $m:=\mathbb{E}[\xi] = rp^{-1}>1$ and thus the process survives forever with positive probability. 

We now choose a random sample of $n$ individuals from the population at some generation $T$, conditional on $n\le Z_T$. A common interest of population genetics is to understand the family tree of the given sample. We can trace back in time the ancestry of such individuals, the lineages will merge in certain generations, until finally all individuals in the sample trace back to a common ancestor. The process that describes such genealogy is known in the literature as the coalescence process and has attracted significant attention in the last decade (see, e.g., \cite{MR4133376, MR4003147, LAMBERT201830, lambertpopovic2013}). 

We consider two types of sampling schemes: Bernoulli sampling and uniform sampling. In the Bernoulli scheme, each individual in generation $T$ is independently selected with a given probability. In contrast, under uniform sampling, a fixed number of individuals is chosen uniformly at random from individuals alive in generation $T$. The purpose of this paper is to study the structure of genealogical trees formed by these kind of samples taken from BGW processes with linear-fractional offspring distributions, and how these genealogical trees are related to each other. Previous work by Lambert and Stadler \cite{lambert2013birth} and Lambert \cite{LAMBERT201830} have analyzed  similar questions in continuous-time settings with binary branching. Here, we shift the focus to discrete time and a linear-fractional offspring reproduction. 
These changes not only introduce new challenges in the proof, but also allows to extend and complete the results of \cite{LAMBERT201830} and \cite{lambert2013birth}.

\section{Coalescent point process}\label{sec:CPP}
We first introduce the coalescent point process associated with $\{Z_n, n\ge 0\}$. In order to do so, we employ a monotone planar embedding of an infinitely old BGW process with arbitrary large population size.   Here, each coordinate $(n,i)$, where $n\in \{0,-1,-2,\dots\}$ and $i\in \{0,1,\dots\}$, corresponds to the $i$-th individual in the $n$-th generation and $\xi_i^{(n)}$ its offspring number. The individual $(0,i)$ will be called the $i$-th tip. We have the following genealogy: the individual $(n,i)$ has mother $(n-1,j)$ if
$$\sum_{k=0}^{j-1}\xi_k^{(n-1)}<i\leq \sum_{k=0}^{j}\xi_k^{(n-1)}.$$
For any $i \ge 1$, denote by $H_i$ the first time when tips $i-1$ and $i$ share a common ancestor, i.e. the coalescent time. The sequence $\{H_i, i\ge 1\}$ forms what is known in literature as the \textit{coalescent point process} (CPP for short). Although the process $\{H_i, i\ge 1\}$ may not be Markovian, in the specific scenario where the offspring distribution is linear-fractional, $\{H_i, i\ge 1\}$ is a sequence of i.i.d. random variables. More precisely, we have the following result. 

\begin{proposition}[Proposition 5.1 in \cite{lambertpopovic2013}]\label{pro:H}
Let $p\in (0,1), r\in [0,1]$ with $p\neq r$. Consider a BGW process $\{Z_n, n\ge 0\}$, whose offspring distribution is $(p,r)$-linear fractional. Then, the branch lengths of the coalescent point process $\{H_i, i\ge 1\}$ are i.i.d. with distribution given by 
\begin{equation}\label{eq:Hdistrib}
 \mathbb{P}(H_i>n)= \frac{r-p}{(1-p)m^n-(1-r)},\qquad n\ge 0,   
\end{equation}
where $m=rp^{-1}$. 
\end{proposition}
In what follows, we work with the CPP $\textbf{H}^{p,r}:=\{H_i, i\ge 1\}$ associated with the BGW  $\{Z_n, n\ge 0\}$ with $(p,r)$-linear fractional offspring distribution.

For any fix $T\in \mathbb{N}$, we denote by $N_T$ the first coalescent time larger than $T$, i.e.
\[N_T:=\min\{i\ge 1: H_i>T\}.\]
Note that   $N_T$ is a geometric random variable with fail probability $\mathbb{P}(H\leq T)$. We will refer to a \textit{coalescent point process with height} $T$ (for short CPP$(T)$) as the sequence $\{H_i, 0\leq i< N_T\}$ where $H_0:=T$.

\begin{figure}
    \centering
        \includegraphics[width=.5\textwidth]{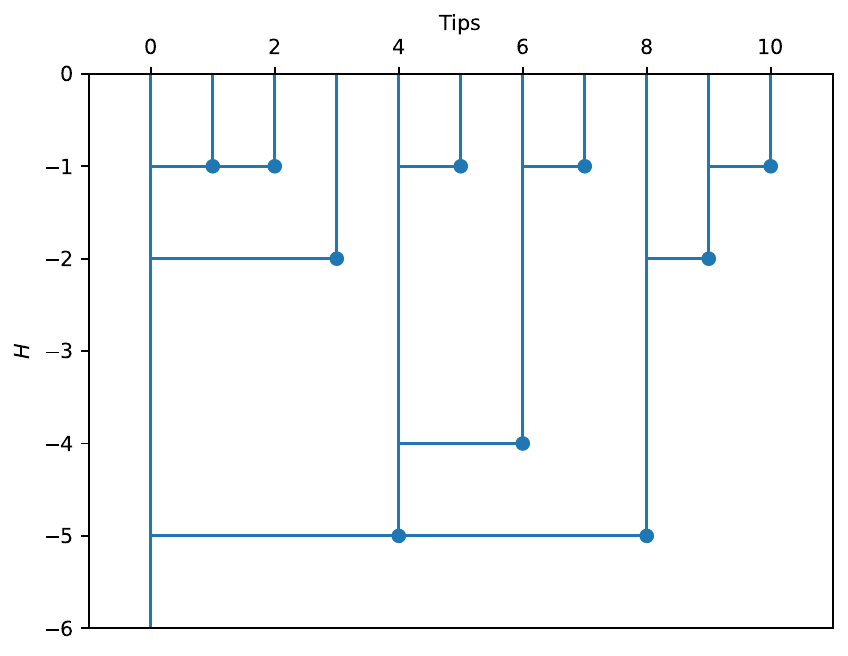}
    \caption{A CPP($6$) tree associated with a $(0.5, 0.8)$-linear fractional BGW process.} 
    \label{fig:coalescents}
\end{figure}

  We can associate a (reduced) tree with the process stopped at the first coalescent time larger than $T$ as follows: initially, we draw a vertical line from $(0,0)$ to $(0,-T)$. Subsequently, we draw additional vertical lines $(i,0)$ to $(i,-H_i)$ for any $i< N_T$. Next, we draw at the bottom of each of these vertical lines, a horizontal line going to the left, stopping when it finds the first vertical line. Note that with this procedure we obtain a ultrametric and oriented tree, with labels from left to right and height in the non-negative integer numbers. 
  We will denote by $\textbf{T}^{p,r}$ the tree constructed as before, and call it a CPP$(T)$ tree. See Figure \ref{fig:coalescents} for simulations of a CPP$(T)$ tree $\textbf{T}^{0.5,0.6}$  with height $T=6$.
Sometimes we also denote the CPP$(T)$ as the tree $\textbf{T}^{p,r}$.  In contrast to the continuous-time setting, we can have multiple coalescent events at the same time. 

From the previous construction, we deduce the following result.

\begin{lemma}\label{likelihoodT}
Let $\tau$ be an ultrametric, oriented tree with $n$ tips and node depths $x_1, \dots, x_{n-1}\in \{1,\dots, T\}$. The likelihood function of the parameters $(p,r)$ given the outcome $\{\mathbf{T}^{p,r}=\tau\}$ is
   \begin{equation*} 
     \mathcal{L}(p, r \mid  \tau) := \mathbb{P}(\mathbf{T}^{p,r}=\tau)=\mathbb{P}(H>T)\prod_{i=1}^{n-1} \mathbb{P}(H = x_i).
   \end{equation*}
Moreover, the likelihood of $(p,r)$ conditioned on $\{N_T=n\}$ is given by
      \begin{equation*}
      \mathcal{L}(p, r \mid  \tau, n): = \mathbb{P}(\mathbf{T}^{p,r}=\tau \mid N_T = n)=\prod_{i=1}^{n-1} \frac{\mathbb{P}(H = x_i)}{\mathbb{P}(H\le T)}.
   \end{equation*}
\end{lemma}

\subsection{The Bernoulli sample and $k$-sample}
In real-world scenarios, it is often not possible to obtain data from the entire population, leading the use of samples instead. If we consider the tree $\mathbf{T}^{p,r}$ with node depths $\{H_{i}, 0\leq i< N_T\}$ where $H_0:=T$, two methods for sampling a population are frequently employed: uniform sampling, where $k$ tips are chosen uniformly with $k\le N_T$; and Bernoulli sampling, where each tip is selected with a given probability $y\in (0,1)$. In both cases, we can construct the subtree associated with the selected tips, which we denoted by $\mathbf{T}_k$ and $\mathbf{T}_y$, respectively. We call $\mathbf{T}_k$ the $k$-sample tree and $\mathbf{T}_y$ the Bernoulli$(y)$ sample tree. For $\mathbf{T}_k$, the associated node depths are denoted by $\{H_{k,i}, 0\leq i< k\}$ where $H_{k,0}:=T$ and $H_{k,i}$ denotes the coalescent time between the $(i-1)$-th and the $i$-th sampled tips, for $i\ge 1$. We note that the node depth $H_{k,i}$ is not necessary the depth of the selected $i$-th tip in $\mathbf{T}^{p,r}$. Similarly, for $\mathbf{T}_y$, the associated node depths are denoted by $\{H_{y,i}, 0\leq i< K\}$ where $H_{y,0}:=T$ and $K$ is the number of selected tips. 

In this section, we are going to study the law of these subtrees and, in particular, how they are related. We start with the Bernoulli$(y)$ sampling with $y\in(0,1)$. First we analyze the Coalescent Point Process (with infinite height) $\textbf{H}^{p,r}=\{H_i, i\geq 1\}$ and sample every tip with probability $y$. The following result gives us the law of $\textbf{H}^{p,r}_{y}:=\{H_{y,i},i\geq 1\}$, where $H_{y,i}$ denotes the coalescent time between the $(i-1)$-th and the $i$-th sampled tips.

\begin{lemma}\label{lemma:Hy} 
Let  $p\in (0,1), r\in [0,1]$ with $p\neq r$ and $y\in(0,1)$. 
The process $\textbf{H}^{p,r}_{y}$ is a Coalescent Point Process associated with a BGW process with $(p_y, r_y)$-linear fractional distribution where $p_y:=1-y(1-p)$ and $r_y:=1-y(1-p)+r-p$. Moreover, the tail distribution of $H_{y,i}$, $i\ge 1$, is given by 
\begin{equation}\label{eq:Hydist}
    \mathbb{P}(H_{y,i}>n)= \frac{r-p}{y(1-p)m^n- [y(1-p)-(r-p)]}, \qquad n\ge  0.
\end{equation}
\end{lemma}

\begin{proof}
The proof directly follows from Proposition \ref{pro:H} and uses analogous arguments to those found in Proposition 2 in \cite{lambert2013birth}. Note that each random variable $H_{y,i}$ has the same distribution as
~$\max\{H_j,1\le j \le J\}$ where $J$ is an independent geometric random variable with success probability $y\in (0,1)$ and the variables $\{H_j, j\geq 1\}$ are i.i.d. random variables with distribution given in \eqref{eq:Hdistrib}. Hence, for $n\geq 0$, 
\begin{equation*}
	\begin{split}
    \mathbb{P}(H_{y,i}\le n) &= \sum_{j=1}^\infty \mathbb{P}\left(\max_{1\le i \le J}H_i\le n \mid J=j\right)\mathbb{P}(J=j) = \sum_{j=1}^\infty \mathbb{P}(H_1\le n)^j y (1-y)^{j-1}\\ &= \frac{y \mathbb{P}(H_1\le n)}{1-(1-y) \mathbb{P}(H_1\le n)}.
    \end{split}
    \end{equation*}
Therefore, for any $n$,
\begin{equation*}
    \mathbb{P}(H_{y,i}> n) =  \frac{\mathbb{P}(H_1> n)}{\mathbb{P}(H_1> n)+y(1- \mathbb{P}(H_1> n))}=\frac{r-p}{y(1-p)m^n- [y(1-p)-(r-p)]}.    \end{equation*}
By comparing with  \eqref{eq:Hdistrib} and noting that the variables $\{H_{y,i},i\geq 1\}$ are independent, we see  that it corresponds to the CPP associated to the BGW process with $(p_y,r_y)$-linear fractional offspring distribution, 
where the parameters are $p_y=1-y(1-p)$ and $r_y=r-p+1-y(1-p)$. 
\end{proof}

\begin{figure}
    \centering
        \includegraphics[width=0.5\textwidth]{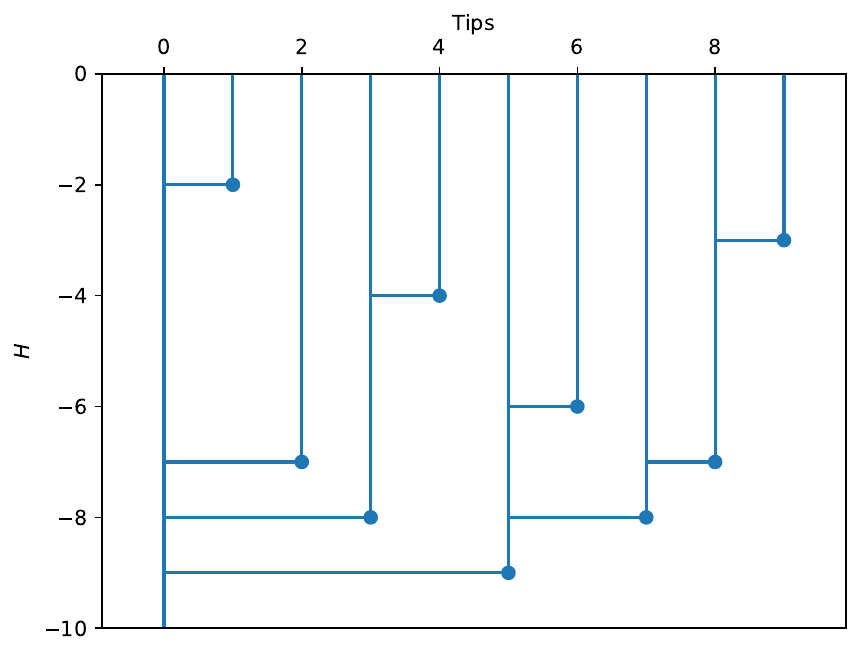}
    \caption{A Bernoulli(0.1) sampled CPP$(10)$ tree  associated with a $(0.5,0.8)$-linear fractional BGW process.} 
    \label{fig:bernoullicoalescents}
\end{figure}

Motivated by the previous result, the process $\textbf{H}^{p,r}_{y}$ will be called a Bernoulli$(y)$ sampled CPP.
In particular, this result tells us that to simulate the genealogy of a Bernoulli$(y)$ sampled tree (infinite) it is sufficient to generate a CPP associated with a $(p_y, r_y)$-linear fractional distribution.  In other words, $\textbf{H}^{p,r}_y= \textbf{H}^{p_y, r_y}$.   
Note that  $\mathbf{T}^{p_y, r_y}$, the reduced associated tree with $\textbf{H}^{p_y, r_y}$,  has the same law as if with probability $y$, we select each tip of the reduced tree $\mathbf{T}^{p,r}$ and construct the associated subtree, i.e. it has the same law as  $\mathbf{T}_y$. See Figure \ref{fig:bernoullicoalescents} for a simulation. 

The next result is an application of Lemma \ref{likelihoodT}, where $K$ is the number of selected tips.

\begin{lemma}\label{lem:likelihoodBernoulli}
Let $\tau$ be an ultrametric, oriented tree with $n$ tips and node depths $x_1, \dots, x_{n-1}\in \{1,\dots, T\}$ and $x_0=T$. The likelihood function of the parameters $(p,r)$ given the outcome $\{\mathbf{T}_y=\tau\}$ is
   \begin{equation*}
    \mathcal{L}(p, r \mid  \tau, y)=
      \mathbb{P}(\mathbf{T}_y=\tau)=\mathbb{P}(H_y>T)\prod_{i=1}^{n-1} \mathbb{P}(H_y = x_i).
   \end{equation*}
Moreover, the likelihood of $(p,r)$ conditioned on $\{K=n\}$ is given by
      \begin{equation*}
      \mathcal{L}(p, r \mid  \tau, y, n):= \mathbb{P}(\mathbf{T}_y=\tau \mid K = n)=\prod_{i=1}^{n-1} \frac{\mathbb{P}(H_y = x_i)}{\mathbb{P}(H_y\le T)}.
   \end{equation*}
\end{lemma}

~\\

Now, we are going to work with the uniform sampling tree $\mathbf{T}_{k}$. Recall that it is constructed as follows.  Consider $\mathbf{T}^{p,r}$ a CPP(T) tree and conditioning on the event $\{N_T>k\}$ with $k\ge 1$, select a uniform sample of $k$ tips. We also recall that the node depths of $\mathbf{T}_{k}$ are denoted by $\{H_{k,i}, 0\leq i< k\}$.

\begin{lemma}\label{lem: density}
Let $m \ge 0$ and $\tau$ be an ultrametric oriented tree with node depths $x_1, \dots, x_{k-1}  \in  \{1,\ldots, T\}$ and  $x_0=T$. The likelihood function of the parameters $(p,r)$ given $\{\mathbf{T}_k = \tau, N_T=k+m\}$ is

\begin{equation*}
\begin{split}
  & \mathcal{L}(p, r\mid \tau,k, m )
  ={ k+m \choose k }^{-1} \underset{\mathbf{m}\in \mathbf{M}_{k,m}}{\sum}\mathbb{P}(H>T)\mathbb{P}(H\leq T)^{m_0} \underset{i=1}{\overset{k-1}{\prod}}
    \left( \mathbb{P}(H\leq x_i)^{m_i+1} -\mathbb{P}(H<x_i)^{m_i+1} \right)  ,
    \end{split}
\end{equation*}
where $\mathbf{M}_{k,m}:=\{(m_0,m_1,\dots, m_{k-1})\in\mathbb{Z}_+^k: m_0+\cdots +m_{k-1}=m\}$.
\end{lemma}

\begin{figure}
    \centering
        \includegraphics[width=\textwidth]{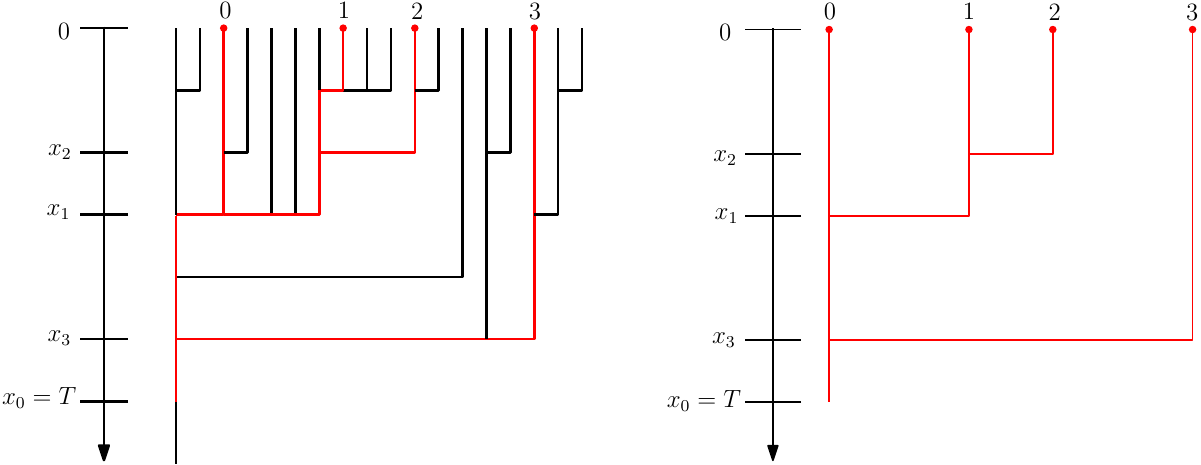}
\caption{Left: A coalescent tree $\mathbf{T}^{p,r}$ of height $T$ and total number of tips $N_T=18$.  The $k$-sampled tree (in red) is included with $k=4$ sampled tips and  corresponding node depths $x_0 = T, x_1, x_2$, and $x_3$. The number of unsampled tips between sampled tips are: $m_2=2$ and $m_0 = m_1 = m_3 = 4$. In each section, the number of tips  with height $x_1,x_2, x_3$  are $n_1 = 3, n_2 = 1$ and $ n_3 = 2$, respectively. Right: The tree $\mathbf{T}_4$ of height $T$ obtained from $\mathbf{T}^{p,r}$ by sampling uniformly 4 tips.
}\label{fig:ksample}
\end{figure}

\begin{proof}
Note that $\mathcal{L}(p, r\mid \tau,k, m )= \mathbb{P}(\mathbf{T}_k=\tau, N_T=k+m)$ is determined by the ratio between favorable cases and possible cases. The total number of ways to select $k$ tips from the $k+m$ possibilities is given by the binomial coefficient ${k+m \choose k}$. Let us denote by $0, 1, 2, \ldots, k-1$ the sampled tips and by $x_i$ the coalescent time between the $(i-1)$-th and the $i$-th sampled tips, for $i\ge 1$.
For each $i\in \{1,\dots, k-1\}$, let $m_i$ be the number of unsampled tips between $i-1$ and $i$, and let $m_0$ represent the number of unsampled tips before $0$ together with the number of unsampled tips after $k-1$. See Figure \ref{fig:ksample}.

Fix an \(i \in \{1, \dots, k-1\}\) and note that the $i$-th sampled tip and all the unsampled tips between \(i-1\) and \(i\) have node depths less than or equal to \(x_i\). Let \(n_i\) denote the number of these tips with depth exactly equal to \(x_i\). More specifically, among the \(m_i+1\) tips between the sampled tip \(i-1\) (excluding it) and the sampled tip \(i\) (including it), there are \(n_i\) tips with depth equal to \(x_i\) and \(m_i+1-n_i\) tips with depths strictly less than \(x_i\), where $1\leq n_i\leq m_i+1$. The term 
\[
\sum_{n_i=1}^{m_i+1} {m_i+1 \choose n_i} \mathbb{P}(H=x_i)^{n_i} \mathbb{P}(H<x_i)^{m_i+1-n_i}
\]  
arises by summing over all possible configurations of the \(n_i\) tips with depth \(x_i\) between the $m_i+1$ tips. From the remains $m_0+1$ tips, the first one has depth strictly bigger than $T$ and the other $m_0$ tips have depth less or equal to $T$. By summing over all the possible vectors $\mathbf{m}\in \mathbf{M}_{k,m}$, the likelihood $\mathcal{L}(p, r\mid \tau,k, m )$ is equal to
\begin{equation*}
\begin{split}
{ k+m \choose k }^{-1} \underset{\mathbf{m}\in \mathbf{M}_{k,m}}{\sum}\mathbb{P}(H>T)\mathbb{P}(H\leq T)^{m_0} \underset{i=1}{\overset{k-1}{\prod}}
    \left(\underset{n_i=1}{\overset{m_i+1}{\sum}} { m_i+1 \choose n_i} \mathbb{P}(H=x_i)^{n_i} \mathbb{P}(H<x_i)^{m_i+1-n_i} \right).
    \end{split}
\end{equation*}

Finally, it is enough to complete the binomial formula for $\mathbb{P}(H\leq x_i)^{m_i+1}$, i.e. to add and subtract the term associated with $n_i=0$ to get the result.
\end{proof}

The previous lemma differs from its continuous-time counterpart  \cite[Section 3]{LAMBERT201830} since in our case we could have multiple branching events at the same time. Sometimes, it is better to have the joint distribution function of the node depths $\{H_{k,i}, 0\le i < k\}$.

\begin{lemma}
Let  $\mathbf{T}_k$ be a $k$-sample tree. For any $x_1,\dots, x_{k-1} \in  \{1,\ldots, T\}$, let $\{y_1,\ldots, y_d\}$ be the set of its different values, $d< k$. Then, for any $m\ge 0$, 
\begin{equation*}
\begin{split}
&\mathbb{P}(N_T=m+k, H_{k,0}=T, H_{k,1}\leq x_1, \dots, H_{k,k-1}\leq x_{k-1})\\
 &\hspace{2.5cm}=
{ k+m \choose k }^{-1} \left( \underset{j=1}{\overset{d}{\prod}}p_j^{r_j}\right)(1-p_0)
 \underset{j=1}{\overset{d}{\sum}} 
 \frac{ p_j^{d-2}(p_j^{m+2}-p_0^{m+2})}{(p_j-p_0)\underset{i=1,i\neq j}{\overset{d}{\prod}}(p_j-p_i)},
\end{split}
\end{equation*}
where $p_0=\mathbb{P}(H\leq T)$, $p_j=\mathbb{P}(H \leq y_j)$ and $r_j=\#\{i <k: x_i=y_j\}$  for $j=1,\ldots,d$.
\end{lemma}

\begin{proof}
We are going to use Lemma \ref{lem: density} and for that we first sum the probability functions up to $x_i$, for every $i=1,\dots, k-1$. More precisely, for each of them, we note that the telescopic sum is equals to 
$$\sum_{y_i=1}^{x_i}\left( \mathbb{P}(H\leq y_i)^{m_i+1} -\mathbb{P}(H<y_i)^{m_i+1} \right)=\sum_{y_i=1}^{x_i}\left( \mathbb{P}(H\leq y_i)^{m_i+1} -\mathbb{P}(H\leq y_i -1)^{m_i+1} \right)=\mathbb{P}(H\leq x_i)^{m_i+1},$$
where we have used that $\mathbb{P}(H\le 0)=0$.
Therefore, by Lemma \ref{lem: density} and the previous computation, we deduce that
\begin{equation*}
\begin{split}
   & \mathbb{P}(N_T=m+k, H_{k,1} \leq x_1, \dots, H_{k,k-1}\leq x_{k-1})\\
    & \hspace{2cm}=
 { k+m \choose k }^{-1} \underset{\mathbf{m}\in \mathbf{M}_{k,m}}{\sum}\mathbb{P}(H>T)\mathbb{P}(H\leq T)^{m_0} \underset{i=1}{\overset{k-1}{\prod}}
    \mathbb{P}(H \leq x_i)^{m_i+1}.
    \end{split}
\end{equation*}
Recall the definition of $r_j$ and observe that
$$\underset{i=1}{\overset{k-1}{\prod}}
   \mathbb{P}(H \leq x_i)^{m_i+1}=\underset{j=1}{\overset{d}{\prod}}
   \mathbb{P}(H \le y_j)^{\widetilde{m}_j+r_j},$$
where $\widetilde{m}_j=\sum_{i: x_i=y_j}m_i$ for $j=1, \dots, d$.
Then 
\begin{equation*}
    \mathbb{P}(N_T=m+k, H_{k,1}\leq x_1, \dots, H_{k,k-1}\leq x_{k-1})=
{ k+m \choose k }^{-1} \sum_{l=0}^m (1-p_0)p_0^{l} 
 \underset{\widetilde{\mathbf{m}}\in \mathbf{M}_{d-1,m-l}}{\sum}\ \underset{j=1}{\overset{d}{\prod}}
  p_j^{\widetilde{m}_j+r_j},
\end{equation*}
where $\widetilde{\mathbf{m}}=(\widetilde{m}_1,\ldots,\widetilde{m}_{d})\in \mathbf{M}_{d-1,m-l}$.
Now, we use Lemma 10 of Lambert and Stadler \cite{lambert2013birth} with 
\[f_{d, m-l}(p_1, \dots, p_d):=\underset{\widetilde{\mathbf{m}}\in \mathbf{M}_{d-1,m-l}}{\sum}\ \underset{j=1}{\overset{d}{\prod}}
  p_j^{\widetilde{m}_j},\] 
  to obtain 
\begin{equation*}
\begin{split}
    \mathbb{P}(N_T=m+k, H_{k,1}\leq x_1, &\dots, H_{k, k-1}\leq x_{k-1})\\
    &=
{ k+m \choose k }^{-1} \left( \underset{j=1}{\overset{d}{\prod}}p_j^{r_j}\right)(1-p_0)
\sum_{l=0}^m  p_0^{l} 
 \underset{j=1}{\overset{d}{\sum}} \frac{ p_j^{m-l+d-1}}{\underset{i=1,i\neq j}{\overset{d}{\prod}}(p_j-p_i)}.
    \end{split}
\end{equation*}
By using the geometric sum, we obtain our result.
\end{proof}

The following result provides a connection between the law of the trees $\mathbf{T}_k$ and $\mathbf{T}_{y}$. Recall that the tree $\mathbf{T}_k$ was constructed from $\{H_i, 0\leq i< N_T\}$, the heights of a CPP$(T)$, by conditioning on $\{N_T\geq k\}$ and uniformly selecting a sample of $k$ tips from $\{0\leq i< N_T\}$.
Define the random variables $\{I_i, 0\leq i< N_T\}$ by
\[
I_i:=\begin{cases}
			1, & \text{if the $i$-th tip is selected,}\\
            0, & \text{otherwise}.
		 \end{cases}
\]
The tree $\mathbf{T}_y$ with $y\in(0,1)$, was constructed from $\{\tilde{H}_i, 0\leq i< \tilde{N}_T\}$, the heights of an independent CPP$(T)$, by selecting independently each tip with probability $y$. Define the random  variables $\{J_i(y), 0\leq i< \tilde{N}_T\}$ by 
\[
J_i(y):=\begin{cases}
			1, & \text{if the $i$-th tip is selected,}\\
            0, & \text{otherwise}.
		 \end{cases}
\]
Observe that $J_i(y)$ is a Bernoulli random variables with parameter $y$. 
Let $K:=\#\{0\leq i< \tilde{N}_T: J_i(y)=1\}$. Denote by $P_y(\cdot)$ the joint law of the random variables 
$\{\tilde{H}_i, 0\leq i< \tilde{N}_T\}, \{J_i(y), 0\leq i< \tilde{N}_T\}$  and $K$.
Let $Y$ be a random variable in $(0,1)$ with improper density $y^{-1}\mathrm{d} y$. Now, we define
\begin{equation}\label{eq:Q}
 \mathbb{Q}(A,Y\in \mathrm{d}y):=y^{-1}P_y(A) \mathrm{d} y,   
\end{equation}
for any  $A$ that can be expressed in terms of  $\{\tilde{H}_i, 0\leq i< \tilde{N}_T\},$ $\{J_i(y), 0\leq i<\tilde{N}_T\}$ and $K$.

The proof of the next lemma follows by using similar arguments as Lemma 6 in \cite{LAMBERT201830}, so we omit the proof.
\begin{lemma}\label{Lem:L}
Let $m\ge k$ and $\textbf{i}=(i_0,\dots, i_{m-1})\in \{0,1\}^m$ where ~$i_0+\dots +i_{m-1}=k$. Then
\begin{equation*}
\begin{split}
        & \mathbb{P}(N_T=m,\,   (I_0,\dots, I_{m-1})= \textbf{i} \mid N_T\ge k)  = \mathbb{Q}(\tilde{N}_T  = m, \, (J_0(Y),\dots, J_{m-1}(Y))= \textbf{i}\mid K=k).
\end{split}
\end{equation*}
Moreover,
\begin{equation}\label{eq:muk}
 \mu_k(\mathrm{d}y):=\mathbb{Q}(Y\in \mathrm{d}y \mid K=k) = \frac{k\delta_T y^{k-1}}{(1-(1-\delta_T)(1-y))^{k+1}}\mathrm{d}y, \qquad y\in (0,1),  
\end{equation}
where 
$\delta_T:=\mathbb{P}(H> T)$.
\end{lemma}

With the above results in hand, we can now show the following result that gives us an explicit formula for computing the law of the tree $\mathbf{T}_k$ in terms of  $\textbf{T}_{Y}$.

\begin{theorem}\label{theo:BernoulliUniform}
Let $k\in \mathbb{N}$ and $y\in (0,1)$. Denote by $\{H_{k,i}, 0\le i < k\}$ and $\{H_{y,i}, 0\leq i< K\}$ the node depths of \textnormal{$\textbf{T}_{k}$} and \textnormal{$\textbf{T}_y$}, respectively. Let $x_0=T$ and $x_1,\dots, x_{k-1}\in\{1,2,\dots, T\}$. Then
$$\mathbb{P}(H_{k,0}=x_0,\dots, H_{k,k-1}=x_{k-1} \mid N_T \ge k)= \int_0^1 \mu_k(\mathrm{d}y) \mathbb{P}(H_{y,0}=x_0, \dots, H_{y,k-1}=x_{k-1}\mid K=k).$$
\end{theorem}

\begin{proof}
Let ~$m \in\mathbb{N}$ with ~$m\geq k$ and ~$\textbf{i}=(i_0, \dots, i_{m-1})\in \{0,1\}^m$ where ~$i_0+\dots +i_{m-1}=k$. In the event ~$\{N_T= m,\, (I_0,\dots, I_{m-1})= \textbf{i}\}$, we have by definition of $\{H_{k,i}, 0\le i < k\}$, that
\[H_{k,0}=\max\{H_0,\dots, H_{l_1}\},\quad  H_{k,1}=\max\{H_{1+l_1},\dots, H_{l_2}\}, \quad \dots\,  \quad H_{k,k-1}=\max\{H_{1+l_{k-1}}, \dots, H_{l_k}\},\]
where ~$l_0:=0$, ~$l_1:= \min\{j\ge 0: i_j\not=0\}$, and ~$l_{r+1}:= \min\{j>l_r: i_j\not=0\}$ for $1\le r < k$.
Similarly, in the event $\{\tilde{N}_T= m,\, (J_0(Y),\dots, J_{m-1}(Y))= \textbf{i}\}$, we have \[H_{Y,0}=\max\{\tilde{H}_0,\dots, \tilde{H}_{l_1}\},\quad  H_{Y,1}=\max\{\tilde{H}_{1+l_1},\dots, \tilde{H}_{l_2}\}, \quad \dots\quad   H_{Y, k-1}=\max\{\tilde{H}_{1+l_{k-1}}, \dots, \tilde{H}_{l_k}\}.\]
Then, by Lemma \ref{Lem:L}, the fact that $H_i\overset{d}{=} \tilde{H}_i$ for any $i\in \mathbb{N}$ and $N_T\overset{d}{=}\tilde{N}_{T}$, we obtain  
 \begin{equation*}
     \begin{split}
      & \mathbb{P}(H_{k,0}=x_0,\dots, H_{k,k-1}=x_{k-1} \mid N_T\ge k) =  \sum_{m\in \mathbb{N}} \sum_{\textbf{i} \in \{0,1\}^m} \mathbb{P}(N_T=m,\,   (I_0,\dots, I_{m-1})= \textbf{i}, \\
      & \hspace{5cm}  \max\{H_0,\dots, H_{l_1}\}=x_0,\ \dots\,  , \max\{H_{1+l_{k-1}}, \dots, H_{l_k}\}=x_{k-1} \mid N_T\ge k)\\ 
      & \hspace{4.6cm} = \sum_{m\in \mathbb{N}} \sum_{\textbf{i}\in \{0,1\}^m} \mathbb{Q}(\tilde{N}_T=m,\,   (J_0(Y),\dots, J_{m-1}(Y))= \textbf{i}, \\
      & \hspace{5cm} \max\{\tilde{H}_0,\dots, \tilde{H}_{l_1}\}=x_0, \ \dots\,  , \max\{\tilde{H}_{1+l_{k-1}}, \dots, \tilde{H}_{l_k}\}=x_{k-1} \mid K=k) \\ & \hspace{4.6cm} = \mathbb{Q} (H_{Y,0}=x_0,\dots, H_{Y, k-1}=x_{k-1} \mid K=k). 
     \end{split}
 \end{equation*}  
 Now, by the total law of probability we obtain 
 \begin{equation*}
     \begin{split}
     & \mathbb{Q}(H_{Y,0}=x_0, \dots, H_{Y,k-1}=x_{k-1}\mid K=k)   \\ & \hspace{1cm} = \int_0^1  \mathbb{Q}(Y\in \mathrm{d}y \mid K=k) \mathbb{Q}(H_{y,1}=x_1, \dots, H_{y,k-1}=x_{k-1}\mid Y=y,  K=k)  \\ & \hspace{1cm} = \int_0^1 \mu_k(\mathrm{d}y) P_y(H_{y,0}=x_0, \dots, H_{y,k-1}=x_{k-1}\mid   K=k),
     \end{split}
 \end{equation*}  
 where the last equality is obtained by definition of the measures $\mathbb{Q}$ and $\mu_k$ given in \eqref{eq:Q} and \eqref{eq:muk}.
\end{proof}

An immediate consequence of the previous result is the following observation, which was stated within the context of a birth and death process in Corollary 4 in \cite{LAMBERT201830}. This result gives a method to simulate the genealogical tree of a $k$-uniform sample through a Bernoulli sampling with probability of success given by the improper distribution of $Y$.

\begin{corollary}\label{lem:consH}
Let $k\in \mathbb{N}$. The law of $H_{k,0},H_{k,1}, \dots, H_{k,k-1}$ can be obtained as follows:
\begin{enumerate}
    \item Select $Y$ with the probability distribution $\mu_k$ defined in \eqref{eq:muk}.
\item Conditioned on $Y=y$, select $H_{y,0}, \dots, H_{y, k-1}$, where $H_{y,0}=T$ and $H_{y,1}, \dots, H_{y, k-1}$ are i.i.d. random variables with distribution 
$$\mathbb{P}(H_{y,i}\leq j\mid H_{y,i}\leq T)= \frac{(m^j-1)}{(m^T-1)}\ \frac{r-p+y(1-p)(m^T-1)}{r-p+y(1-p)(m^j-1)}, \qquad 1\leq j\leq T,\  1\leq i < k,$$
where $m=r/p$.
\end{enumerate}
Furthermore, let $\tau$ be an ultrametric, oriented tree with node depths $x_0=T$ and $x_1, \dots, x_{k-1}  \in  \{1,\ldots, T\}$, the likelihood function of the parameters $(p,r)$ given the outcome $\{\mathbf{T}_k=\tau\}$ is
   \begin{equation*}
    \mathcal{L}(p, r,\mid  \tau,k) = \mathbb{P}(\mathbf{T}_k=\tau)=
       \int_0^1 \mu_k(\mathrm{d}y) \prod_{i=1}^{k-1} \frac{\mathbb{P}(H_{y,i} = x_i)}{\mathbb{P}(H_{y,i}\le T)}.
   \end{equation*}
\end{corollary}

\begin{proof}
The first part holds by Theorem \ref{theo:BernoulliUniform} and Lemma \ref{lemma:Hy}. For the likelihood function, we use Theorem~\ref{theo:BernoulliUniform} and Lemma \ref{lem:likelihoodBernoulli}.
\end{proof}

\section{Discussion}

\textbf{Impact for parameter inference.} 
In this section, we illustrate how our results can be applied in practice. Although real data were not available to us, we would like to point it out how our results can be used for inference purposes on real-world datasets. 

The classical maximum likelihood approach can be used to estimate the model parameters $r$ and $p$, which characterize the offspring distribution of the BGW process. Suppose we observe a dataset of the form
\[x_1,\dots, x_{n-1}\in \{1, \dots, T\},\]
where the values $x_1, \dots, x_{n-1}$ represent an i.i.d. sample drawn from the distribution of the coalescent time $H$, conditioned on $\{H \leq T\}$. In other words, these are the coalescent times observed in a BGW tree with a linear fractional offspring distribution with certain parameters $r$ and $p$.

With a set of data like this and Lemma \ref{likelihoodT}, the true parameters of the whole tree can be estimated. Under the CPP($T$) model, the likelihood of observing a particular tree shape with these coalescent depths is $ \mathcal{L}(p, r \mid  \tau)$ given in Lemma \ref{likelihoodT}. Using the maximum likelihood, the parameters of the model $r$ and $p$ can be estimated numerically. Specifically, the maximum likelihood estimates $(\hat{p}, \hat{r})$ are obtained as
\[(\hat{p}, \hat{r}) = \arg\max_{(p,r)}  \mathcal{L}(p, r \mid  \tau),\]
where, by using Proposition \ref{pro:H} and  Lemma \ref{likelihoodT}, the likelihood is given by
\[ \mathcal{L}(p, r \mid  \tau)=\frac{p^T(r-p)^n}{(1-p)r^T-(1-r)p^T} \prod_{i=1}^{n-1} \left(\frac{p^{x_i-1}}{(1-p)r^{x_i-1}-(1-r)p^{x_i-1}}- \frac{p^{x_i}}{(1-p)r^{x_i}-(1-r)p^{x_i}}\right).\]
Since the maximum likelihood estimator does not admit a closed-form expression, it must be computed using numerical optimization methods.

In practice, datasets of the form described above, i.e. the coalescent times of all individuals in the population, may not be directly available. More commonly, one may only have access to a sample taken at a fixed generation under a Bernoulli or an uniform sampling scheme. In such cases, inference of the parameters can still be performed using the likelihood expressions given in Lemma \ref{lem:likelihoodBernoulli} and Corollary \ref{lem:consH} for Bernoulli and uniform sampling, respectively. 
\\ 

\textbf{Connection with the CPP for birth-death processes.} We can embed the BGW process $\{Z_n, n\ge 0\}$ with $(p,r)$-linear fractional offspring distribution in a birth and death process.  By using Chapter III, Section 6 in \cite{MR0373040}, we have the following simple observation.

\begin{lemma}\label{Theo:embedding}
    Let $\{Z_n, n\ge 0\}$ be a BGW process with $(p,r)$-linear fractional offspring distribution, $p\neq q$. Then $\{Z_n, n\ge 0\}$ is embedded in a birth and death process $\{Z_t, t\ge 0\}$ initiated with a single individual and with birth and death rates given by 
    $$\lambda=(1-p)\frac{\log(p)-\log(r)}{p-r} \qquad \mbox{and} \qquad  \mu=(1-r)\frac{\log(p)-\log(r)}{p-r}.$$
\end{lemma}

This embedding helps clarify the analogy between our results and those established in \cite{LAMBERT201830}, where a birth–death process in  continuous-time is considered to construct genealogical trees. In this setting, there is only one branching event at a given time. In contrast, our setting differs in both time and offspring structure: multiple individuals may reproduce simultaneously and potentially generating more than two offspring. This leads to different genealogical trees.

In applications involving real-world data, such as the  periodical cicadas, a discrete-time framework seems to provide a more natural modeling choice. Rather than embedding the observed discrete tree structure into a continuous-time process, it may be more faithful to the biological reality to analyze the genealogies directly through a discrete-time model. This approach allows to  complement the results from \cite{LAMBERT201830}, extending the study of genealogies beyond the binary tree assumption and the continuous-time models.

\section*{Acknowledgments}
S. Palau acknowledges support from UNAM-PAPIIT (IN103924).

\end{document}